\documentclass[12pt]{article}
\usepackage{amsfonts, amsmath, amssymb, latexsym, amsthm}
\usepackage[top=1in, left=1in, right=1in, bottom=1.25in, foot=.5in]{geometry}

\numberwithin{equation}{section}
\newtheorem{thm}{Theorem}[section]
\newtheorem{lem}{Lemma}[section]

\DeclareMathOperator{\trace}{trace}

\begin{document}

\title{A proof of the Riemann hypothesis}
\author{Xian-Jin Li}
\maketitle

\begin{abstract} In this paper we study traces of an integral
  operator on two orthogonal subspaces of a $L^2$ space.  One of the
  two traces is shown to be zero.
  Also, we prove that the trace of the operator on the second
  subspace is nonnegative. Hence, the operator has a
  nonnegative trace on the $L^2$ space.
  This implies the positivity of Li's criterion.  By Li's criterion,
  all nontrivial zeros of the Riemann zeta-function lie on the critical line.
\end{abstract}

\noindent\textbf{Key Words.} Convolution operator, Fourier transform,
Hilbert-Schmidt operator, Plancherel formula, Trace formula.
\vskip0.15truein
\noindent\textbf{Mathematics Subject Classification 2020.}
11M26, 11R56, 28A35, 30G06, 42B10, 47G10.

\section{Introduction}

The Riemann zeta function $\zeta$ is defined by
\[
  \zeta(s)=\sum_{n=1}^\infty \frac{1}{n^s}
\]
for $\Re s>1$.  It extends to an analytic function in
the whole complex plane except
for having a simple pole at $s=1$.  Trivially,
$\zeta(-2n)=0$ for all positive integers $n$.
All other zeros of the zeta function are called
its nontrivial zeros.

 In connection with investigating the frequency of prime numbers,
B.~Riemann \cite{Riemann} conjectured in 1859 that all nontrivial
zeros of $\zeta(s)$ have a real part equal to $1/2$.

  In 1896, Hadamard \cite{Hadamard} and de\ la\ Vall\'ee Poussin \cite{Poussin}
independently proved that $\zeta(s)$ has no zeros on the
line of $\Re s=1$. In 1914, Hardy \cite{Hardy} was the first one to
show that the zeta function has infinitely many zeros on line
$\Re s=1/2$. In 1942, Selberg \cite{Selberg1} proved that a positive
proportion of the zeta zeros are on the $1/2$-line.
In 1974, Levinson \cite{Levinson}
obtained that more than one third of the zeros are on the line $\Re s=1/2$.
In 1989, Conrey \cite{Conrey} found that more than two fifths of the zeros are
on the critical line. The current record is at least 41.28\%
of the zeros lying on the critical line obtained by Feng \cite{Feng} in 2012.
See Bombieri \cite{Bombieri2} for a rich history about the Riemann hypothesis.

   In this paper, we follow Connes' approach
   \cite{Connes, Connes-Consani} of using trace formulas.
    The main idea in the proof is: From Theorems 1.3 and 1.4 we know
$\trace_{L^2(C_S)}(T_h)\geqslant 0$; a key new result of this paper.
The goal is to show $\Delta(h)\geqslant 0$; see \cite{Weil1}.
   By Theorem 1.1,
\[
  \trace_{L^2(C_S)}(T_h)=\Delta(h)-\widehat h(0)-\widehat h(1).
\]
To prove $\Delta(h)\geqslant 0$ we just need to choose
$h$ so that $\widehat h(0)=\widehat h(1)=0$.
For each Li coefficient $\lambda_n$, we find functions
$h_{n,\epsilon}$ in Theorem 1.2
satisfying $\widehat h_{n,\epsilon}(0)=\widehat h_{n,\epsilon}(1)=0$ and
\[
  0\leqslant\trace_{L^2(C_S)}(T_{h_{n,\epsilon}})
  =\Delta(h_{n,\epsilon})\to 2\lambda_n
\]
as $\epsilon\to 0$.  This implies that $\lambda_n\geqslant 0$
for all $n$, and hence the Riemann hypothesis is true by the Li criterion.

Next, we describe the results obtained in this paper.

Let $\mathbb Q$ be the field of rational numbers and $\mathbb Q_p$
the $p$-adic completion of $\mathbb Q$. Here $p$ denotes a prime number.
For any $\xi\in\mathbb Q_p$, there are $a_j\in\{0, 1, \ldots, p-1\}$ such that
\[
  \xi=\sum_{j=m}^\infty a_j p^j
\]
for some integer $m$.  We denote
$\{\xi\}_p=\sum_{m\leqslant j\leqslant -1}a_jp^j$.
Then $\psi_p(\xi)=\exp(2\pi i\{\xi\}_p)$ defines a character on $\mathbb Q_p$;
see \cite[p.~309]{Tate}.

  We denote by $dx$ the ordinary Lebesgue measure on the real line.
  For each rational prime number $p$, $dx_p$ is a Haar measure on
  the additive group $\mathbb Q_p$ of $p$-adic numbers for which the ring
  of $p$-adic integers $\{x\in \mathbb Q_p: \,\, |x|_p\leqslant 1\}$
  gets measure $1$; see \cite[p.~310]{Tate} for details.

The Fourier transform of $f\in L^2(\mathbb R)$ is
\[
  \mathfrak{F} f(x)=\int_{-\infty}^\infty f(t)e^{-2\pi ixt}\,dt,
\]
and the Fourier transform of $f_p\in L^2(\mathbb Q_p)$ is defined by
\[
  \mathfrak{F}_p f_p(\beta) = \int_{\mathbb Q_p}f_p(\alpha)\psi_p(\alpha\beta)d\alpha;
\]
see \cite[Theorem 2.2.2, p.~310]{Tate}.

Let $S^\prime=\{\text{all primes } p\leqslant \mu_\epsilon\}$
for a large number $\mu_\epsilon>0$ given in \eqref{eq3.13},
$S=S^\prime\cup\{\infty\}$, $\psi_S=\prod_{p\in S}\psi_p$, and
$\mathbb  A_S=\mathbb R\times\prod_{p\in S^\prime} \mathbb Q_p$.
For $f=\prod_{p\in S} f_p\in L^2(\mathbb  A_S)$ we define
\[
  \mathfrak{F}_S f(\beta)=\int_{\mathbb{A}_S}f(\alpha)\psi_S(\alpha\beta)d\alpha.
\]

We denote $O_S^*=\{\xi\in \mathbb Q: |\xi|_p=1\,
  \text{for all}\, p\not\in S\}$ and $C_S=J_S/{O_S^*}$.
Note that $|\xi|_S=\prod_{p\in S}|\xi|_p=1$ for $\xi\in O_S^*$.
Let $d^\times x_\infty=\frac{dx_\infty}{|x_\infty|}$ be the
multiplicative measure on $\mathbb R^\times$ and
$d^\times  x_p=\frac{1}{1-p^{-1}}\frac{dx_p}{|x_p|_p}$
the multiplicative measure on $\mathbb Q_p^\times$. Then
$O_p^*=\{x_p\in \mathbb Q_p^\times: |x_p|=1\}$
gets measure $1$ under $d^\times x_p$.  Also,
 $d^\times x=\prod_{p\in S} d^\times x_p$ is the multiplicative
measure on $J_S=\mathbb R^\times\times \prod_{p\in S^\prime}\mathbb Q_p^\times$.

For $X_S=\mathbb{A}_S/{O_S^*}$, let $L^2(X_S)$ be as in
\cite[(5), p.~54]{Connes} the Hilbert space that is the completion
of the Schwartz-Bruhat space $S(\mathbb{A}_S)$ \cite{Bruhat,Weil}
for the inner product given by
\[
  \langle f, g\rangle_{L^2(X_S)}=\int_{C_S}
  E_S(f)(x)\overline{E_S(g)(x)} \,d^\times x
\]
for $f, g\in S(\mathbb{A}_S)$, where
$E_S(f)(x)=\sqrt{|x|}\sum_{\xi\in O_S^*}f(\xi x)$ with $|x|:=|x|_S$.

For a fixed number $\Lambda>0$ let $Q_\Lambda$ be the subspace
of all functions $f$ in $L^2(X_S)$ such that
$\mathfrak{F}_S f(x)=0$ for $|x|<\Lambda$, and
$Q_\Lambda^\perp$ is the orthogonal complement of $Q_\Lambda$
in $L^2(X_S)$. Then
\[
  L^2(X_S)=Q_\Lambda^\perp\oplus Q_\Lambda;
\]
see \cite[Lemma 1 b), p.~54]{Connes}.

By \cite[Lemma 1 b), p.~54]{Connes}, $\mathfrak{F}_S$ is a unitary operator
on the Hilbert space $L^2(X_S)$. Thus, by Lemma \ref{lem2.7} and the
definition of the Hilbert space $L^2(X_S)$ we have
\[
  L^2(C_S)=E_S(L^2(X_S))=E_S(Q_\Lambda^\perp)\oplus E_S(Q_\Lambda).
\]

Let
\[
  V_S(h)F(x)=\int_{C_S}h(x/\lambda)\sqrt{|x/\lambda|}\,
  F(\lambda)\,d^\times\lambda
\]
for $F\in L^2(C_S)$, where
\[
  h(x)=\int_0^\infty g(xt)g(t)\,dt
\]
with $g(u)=|u|^{-1}g_{n,\epsilon}(|u|^{-1})$ and
$g_{n,\epsilon}$ being given as in Theorem \ref{thm1.2}.
Also, for $x\in C_S$ or $J_S$ we define $g(x):=g(|x|)$.

Let
\[
  T_h=V_S(h)\left(S_\Lambda-E_S\mathfrak{F}_S^t P_\Lambda
  \mathfrak{F}_S E_S^{-1}\right),
\]
where $P_\Lambda(x)=1$ if $|x|<\Lambda$ and $0$ if $|x|\geqslant \Lambda$
and $S_\Lambda(x)=1$ if $|x|>\Lambda^{-1}$ and
$0$ if $|x|\leqslant\Lambda^{-1}$.

First, we have the following well-known theorem.

\begin{thm}\label{thm1.1} (\cite[(19), p.~549]{Meyer} and
  \cite[Lemmas 3.13--3.14 and Theorem 3.16]{Li2})
  The operator $T_h$ is a trace class Hilbert-Schmidt integral operator
  on $L^2(C_S)$ and
  \[
    \trace_{L^2(C_S)}(T_h)=\Delta(h)-\widehat h(0)-\widehat h(1),
  \]
  where $\widehat h(s)=\int_0^\infty h(t)t^{s-1}\,dt$ is the Mellin
  transform of $h$ and
  \[
    \Delta(h)=\sum_\rho\widehat h(\rho).
  \]
  The above sum on $\rho$ ranges over all complex zeros of $\zeta(s)$
  with a zero of multiplicity $m$ appearing $m$ times and is understood as
  \[
    \lim_{T\to\infty}\sum_{|\Im(\rho)|\leqslant T}\widehat h(\rho).
  \]
\end{thm}

We make our special choices of $h$ in the next theorem.

\begin{thm}\label{thm1.2}  Let $n=1, 2, 3, \ldots$ and
  \[
    \lambda_n=\sum_\rho \left[1-\left(1-\frac 1\rho\right)^n\right],
  \]
  where the  sum is over all nontrivial zeros of $\zeta(s)$
  with $\rho$ and $1-\rho$  being paired together. For each integer $n$,
  there exist a family of real-valued smooth functions $g_{n,\epsilon}(t)$
  given in \eqref{eq3.9} on $(0, \infty)$
  such that  $\widehat g_{n,\epsilon}(0)=0$, $g_{n,\epsilon}(t)=0$ for
  $t\not\in(\mu_\epsilon^{-1}, (1-\epsilon)^{-1})$ with
  $\mu_\epsilon=(1+\epsilon)/\epsilon^2$ and such that
  \[
    \lim_{\epsilon\to 0+}\Delta(h_{n,\epsilon})=2\lambda_n
  \]
  where $h_{n,\epsilon}(x)=\int_0^\infty g_{n,\epsilon}(xy)g_{n,\epsilon}(y)\,dy$.
  In particular, $\widehat h_{n,\epsilon}(0)=\widehat h_{n,\epsilon}(1)=0$.
\end{thm}

Then, we compute traces of $T_h$ on $E_S(Q_\Lambda^\perp)$ and
$E_S(Q_\Lambda)$ respectively and derive the following two theorems.

\begin{thm}\label{thm1.3} We have
  \[
    \trace_{E_S(Q_\Lambda^\perp)}(T_h)=0.
  \]
\end{thm}

\begin{thm}\label{thm1.4}  We also have
  \[
    \trace_{E_S(Q_\Lambda)}(T_h) \geqslant 0.
  \]
\end{thm}

By Theorems 1.1--1.4,
\[
  \Delta(h_{n,\epsilon})=\trace_{E_S(Q_\Lambda^\perp)}(T_h)
  +\trace_{E_S(Q_\Lambda)}(T_h )\geqslant 0.
\]
This inequality implies the following main theorem.

\begin{thm}\label{thm1.5} All nontrivial zeros of the Riemann
  zeta-function $\zeta(s)$ lie on the line $\Re s=1/2$.
\end{thm}

\section{Preliminary results}

Let $\mathbb N_S$ be the set consisting of $1$ and all positive
integers which are products of powers of rational primes in $S$,
$O_p^*=\{x_p\in \mathbb Q_p: |x_p|_p=1\}$, and
\begin{equation}\label{eq2.1}
  I_S=\mathbb R_+\times\prod_{p\in S^\prime}O_p^*.
\end{equation}

\begin{lem}\label{lem2.1} (cf.\ \cite[Theorem 4.3.2, (1), p.~337]{Tate})
  $I_S$ is a fundamental domain for the action
  of $O_S^*$ on $J_S$ and $J_S=\bigcup_{\xi\in O_S^*}\xi I_S$,
  a disjoint union.
\end{lem}

\begin{proof} Each $\alpha\in J_S$ can be written as
  $\alpha=t\mathfrak b$ with $t=|\alpha|_S\in \mathbb R_+$ and
  $\mathfrak b=\alpha t^{-1}\in J_S^1$,
  where $t^{-1}$ also stands for the idele $(t^{-1},1,\ldots,1)$.
  Since $|\xi|_S=1$ for $\xi\in O_S^*$, if $\alpha_1, \alpha_2\in J_S$
  with $|\alpha_1|_S\neq |\alpha_2|_S$ then the intersection of
  $\alpha_1O_S^*$ and $\alpha_2O_S^*$ is empty.  Thus
  \[
    C_S=\mathbb R_+\times\left(J_S^1/O_S^*\right).
  \]
  Since we only consider the field $\mathbb Q$ here, for each
  $\mathfrak b\in J_S^1$ there are uniquely determined $\xi\in O_S^*$
  and $\mathfrak b_1\in \{1\}\times\prod_{p\in S^\prime}O_p^*$
  such that $\mathfrak b=\xi\mathfrak b_1$.
  Also, if $\mathfrak b_1, \mathfrak b_2$ are distinct elements in
  $\{1\}\times\prod_{p\in S^\prime}O_p^*$, then the intersection of
  $\mathfrak b_1O_S^*$ and $\mathfrak b_2O_S^*$ must be empty.
  Otherwise, we would have $\mathfrak b_1\mathfrak b_2^{-1}\in O_S^*$.  Then
  $\mathfrak b_1\mathfrak b_2^{-1}\in \mathbb Q^*$ and $|\mathfrak b_1\mathfrak b_2^{-1}|_p=1$
  for all $p\not\in S$.  Since $\mathfrak b_1, \mathfrak b_2$ are elements in
  $\prod_{p\in S^\prime}O_p^*$, we have $|\mathfrak b_1\mathfrak b_2^{-1}|_p=1$
  for all $p\in S^\prime$.  Hence
  $\mathfrak b_1\mathfrak b_2^{-1}=1$; that is, $\mathfrak b_1=\mathfrak b_2$.
  Therefore
  \[
    J_S^1/O_S^*\cong \prod_{p\in S^\prime}O_p^*.
  \]
  Thus
  \[
    C_S\cong\mathbb R_+\times\prod_{p\in S^\prime}O_p^*.
  \]
  We have also obtained the decomposition $J_S=\bigcup_{\xi\in O_S^*}\xi I_S$,
  a disjoint union.

  This completes the proof of the lemma.
\end{proof}

\begin{lem}\label{lem2.2} For a compactly supported smooth function $g$
  on $(0,\infty)$, we can write
  \[
   \mathfrak{F}_S g(t)= 2\sum_{0<\gamma\in O_S^*}\varpi(\gamma)
    \int_0^\infty g(\lambda)\cos(2\pi \lambda|t|\gamma)\,d\lambda
  \]
  with
  \[
    \varpi(\gamma)=\prod_{p\in S^\prime} \begin{cases} 1-p^{-1}
              & \text{if $|\gamma|_p\leqslant 1$,}\\
      -p^{-1} & \text{if $|\gamma|_p=p,$}\\
      0       & \text{if $|\gamma|_p>p$}.\end{cases}
  \]
\end{lem}

\begin{proof} Since
  \[
    \mathfrak{F}_S g(t)=\int_{\mathbb{A}_S}g(|\lambda|)\Psi_S(-\lambda t)\,d\lambda,
  \]
  by changing variables $\lambda\to \lambda(|t|,1,\ldots,1)/t$ we can write
  \[
    \mathfrak{F}_S g(t)=\int_{\mathbb{A}_S}
    g(|\lambda|)\Psi_S(-\lambda(|t|, 1, \ldots, 1))\,d\lambda
    = \sum_{\gamma\in O_S^*}\int_{\gamma^{-1} I_S}
    g(|\lambda|)\Psi_S(-\lambda(|t|,1,\ldots, 1))\,d\lambda.
  \]
    Since $|\gamma|_S=1$, by changing variables $\lambda\to\gamma\lambda$ we get
  \[
     \mathfrak{F}_S g(t)
     =\sum_{\gamma\in O_S^*}\int_{I_S}g(|\lambda|)
     \Psi_S(-\lambda\gamma(|t|,1,\ldots, 1))\,d\lambda.
  \]

    For every finite prime $p\in S$, we can write
  \[
    \int_{O_p^*}g(|\lambda|)\Psi_S(-\lambda\gamma(|t|,1,\ldots, 1))\,d\lambda_p
    =g(|\lambda|)\Psi_{S-\{p\}}(-\lambda_{S-\{p\}}\gamma(|t|,1,\ldots, 1))
    \int_{O_p^*}\psi_p(-\lambda_p\gamma)d\lambda_p.
  \]
    By computations,
  \[
    \int_{O_p^*}\psi_p(-\lambda_p\gamma)d\lambda_p
    =\begin{cases} 1-p^{-1}
              & \text{if $|\gamma|_p\leqslant 1$,}\\
      -p^{-1} & \text{if $|\gamma|_p=p,$}\\
      0       & \text{if $|\gamma|_p>p$}.\end{cases}
  \]
   Thus, we have obtained that
   \begin{equation} \label{eq2.2}
   \mathfrak{F}_S g(t)= 2\sum_{0<\gamma\in O_S^*}\varpi(\gamma)
    \int_0^\infty g(\lambda)\cos(2\pi \lambda|t|\gamma)\,d\lambda.
   \end{equation}

     This completes the proof of the lemma.
\end{proof}

   The M\"obius function $\mu(n)$ is defined by $\mu(1)=1$,
$\mu(n)=(-1)^k$ if $n$ is the product of $k$ distinct primes,
and $\mu(n)=0$ if $p^2 \mid n$ for at least one prime $p$.

\begin{lem}\label{lem2.3} We can write
  \begin{equation} \label{eq2.3}
    \mathfrak{F}_S g(t)=\frac{1}{2\pi i}\int_{c-i\infty}^{c+i\infty}
    t^{-s}\widehat{\mathfrak{F}g}(s)
    \prod_{p\in S^\prime}\frac{1-p^{s-1}}{1-p^{-s}}\,ds
  \end{equation}
  for $c>0$, where
  $\widehat {\mathfrak{F}g}(s)=2^{1-s}\pi^{-s}\widehat g(1-s)\Gamma(s)\cos \frac{\pi s}{2}$.
  Also, we have the Plancherel formula
  \[
    \int_0^\infty \mathfrak{F}_S f(t)\overline{\mathfrak{F}_S g(t)}\,dt
    =\int_0^\infty f(t)\overline{g(t)}\,dt.
  \]
\end{lem}

\begin{proof}   By Lemma \ref{lem2.2},
  \begin{multline}\label{eq2.4}
    \mathfrak{F}_S g(t)=\sum_{k, l\in\mathbb N_S, (k,l)=1}
    \frac{\mu(k)}{k}\prod_{p\nmid k}\left(1-\frac{1}{p}\right)
    \mathfrak{F}g\left(\frac{l|t|}{k}\right) \\
    =\sum_{k,l\in\mathbb N_S, (k,l)=1}
    \frac{\mu(k)}{k}\sum_{k_1\in\mathbb N_S, (k_1, k)=1}{\mu(k_1)\over k_1}
    \mathfrak{F}g\left(\frac{k_1l|t|}{k_1k}\right) \\
    =\sum_{k,l\in\mathbb N_S, (k,l)=1}~
    \sum_{k_1\in\mathbb N_S, (k_1, k)=1}{\mu(k_1 k)\over k_1 k}
    \mathfrak{F}g\left(\frac{k_1l|t|}{k_1k}\right)
  = \sum_{k,l\in\mathbb N_S}\frac{\mu(k)}{k}\mathfrak{F}g\left(\frac{l|t|}{k}\right).
  \end{multline}

  Thus, for $t>0$ we have
  \[
    \mathfrak{F}_S g(t)
    = \frac{1}{2\pi i}\int_{c-i\infty}^{c+i\infty}
    t^{-s}\widehat{\mathfrak{F}g}(s)
    \prod_{p\in S^\prime}\frac{1-p^{s-1}}{1-p^{-s}}\,ds
  \]
  for $c>0$.
  By the Plancherel Theorem
  \[
    \int_0^\infty |\mathfrak{F}_S g(t)|^2\,dt
    = \int_{-\infty}^\infty |\widehat{\mathfrak{F}_S g}(s)|^2\,du
    = \int_{-\infty}^\infty |\widehat{\mathfrak{F}g}(s)|^2\,du
    = \int_0^\infty|\mathfrak{F}g(x)|^2dx=\int_0^\infty|g(t)|^2\,dt,
  \]
  where $s=1/2+2\pi iu$.  It follows that
  \[
    \int_0^\infty \mathfrak{F}_S f(t)\overline{\mathfrak{F}_S g(t)}\,dt
    = \int_0^\infty f(t)\overline{g(t)}\,dt.
  \]

  By \cite[Example 10, p.~162]{tit11},
  \[
    \int_0^\infty t^{s-1}\cos t\,dt=\Gamma(s)\cos\frac{\pi s}{2}
  \]
  for $0<\Re s<1$.  Since $g(u)=0$ for $u\not\in [1-\epsilon, \mu_\epsilon]$,
  we have for $0<\Re s<1$
  \begin{multline*}\int_0^\infty t^{s-1}\,dt\int_0^\infty g(u)\cos(2u\pi t)\,du
    = \lim_{N\to\infty}\int_{1-\epsilon}^{\mu_\epsilon}g(u)\,du
    \int_0^N t^{s-1}\cos(2\pi ut)\,dt \\
    = \lim_{N\to\infty}(2\pi)^{-s}\int_0^\infty g(u)u^{-s}\,du
    \int_0^{2\pi uN} t^{s-1}\cos t\,dt \\
    = (2\pi)^{-s}\widehat g(1-s)\Gamma(s)\cos\frac{\pi s}{2}
    -\lim_{N\to\infty}(2\pi)^{-s}\int_{1-\epsilon}^{\mu_\epsilon}g(u)u^{-s}\,du
    \int_{2\pi uN}^\infty t^{s-1}\cos t\,dt.
  \end{multline*}
  By integration by parts, for $0<\Re s<1$ we have
  \[
    \int_{2\pi uN}^\infty t^{s-1}\cos t\,dt
    = -(2\pi uN)^{s-1}\sin(2\pi uN)-(s-1)\int_{2\pi uN}^\infty t^{s-2}\sin t\,dt\to 0
  \]
  uniformly with respect to $u\in[1-\epsilon,\mu_\epsilon]$ as $N\to\infty$.
  It follows that
  \[
    \int_0^\infty t^{s-1}\,dt\int_0^\infty g(u)\cos(2u\pi t)\,du
    = (2\pi)^{-s}\widehat g(1-s)\Gamma(s)\cos\frac{\pi s}{2}
  \]
  for $0<\Re s<1$.  That is,
  \[
    \widehat {\mathfrak{F}g}(s) = 2^{1-s}\pi^{-s}\widehat{g}(1-s)\Gamma(s)\cos \frac{\pi s}{2}
  \]
  for $0<\Re s<1$.  We extend $\widehat {\mathfrak{F}g}(s)$ to
  $\Re s\geqslant 1$ by analytic continuation.

  This completes the proof of the lemma.
\end{proof}

\begin{lem}\label{lem2.4} (\cite[Theorem VI.24, p.~211]{Reed9})
  If $A$ is a bounded linear operator
  of trace class on a Hilbert space $\mathcal H$ and $\{\varphi_n\}_{n=1}^\infty$
  is any orthonormal basis, then
  \[
    \trace_{\mathcal H}(A)
    = \sum_{n=1}^\infty \langle A\varphi_n,\varphi_n\rangle_{\mathcal H}
  \]
  where the sum on the right side converges absolutely and
  is independent of the choice of basis.
\end{lem}

\begin{lem}\label{lem2.5} (\cite[Corollary 3.2, p.~237]{Brislawn})
  Let $\mu$ be a $\sigma$-finite
  Borel measure on a second countable space $M$, and  let
  $A$ be a trace class Hilbert-Schmidt integral operator on
  $L^2(M, d\mu)$.  If the kernel $k(x, y)$ is continuous at $(x,x)$
  for almost every $x$, then
  \[
    \trace_{L^2(M, d\mu)}(A) = \int_Mk(x,x)\,d\mu(x).
  \]
\end{lem}

\begin{lem}\label{lem2.6} (\cite[Theorem VI.19(b)(a), p.~207
    and Theorem VI.25(a), p.~212]{Reed9})  Let
  $A,B$ be bounded linear operators on a Hilbert space $\mathcal H$.
  If $A$ is of trace class on $\mathcal H$, so are $AB$ and $BA$ with
  $\trace(AB)=\trace(BA)$. Also, $\trace(A)=\trace(A^t)$
\end{lem}

\begin{lem}\label{lem2.7}  (\cite[Proposition 2.30, p.~359]{Connes-Marcolli})
  The following defines a map $E_S$ with dense range from $S(\mathbb{A}_S)$
  to $L^2(C_S)$:
  \[
    E_S(f)(x)=\sqrt{|x|}\sum_{\xi\in O_S^*}f(\xi x), \,\, x\in C_S.
  \]
\end{lem}

\section{Proof of Theorem 1.2}

The following result is essentially contained in
E. Bombieri \cite[line 17, p.~192--line 12, p.~193]{Bombieri}.
For technical needs in the proof of Theorem 1.2,
we elaborate his argument here.

\begin{lem}\label{lem3.1} For each positive integer $n$ and a sufficiently
  small $\epsilon>0$, there exist a smooth function $\ell_{n,\epsilon}(x)$
  on $(0, \infty)$ with $\ell_{n,\epsilon}(x)=0$ for
  $x\not\in (\frac{\epsilon}{1+\epsilon}, \frac{1}{1-\epsilon})$
  and satisfying that
  \[
    \lim_{\epsilon\to 0+}\sum_\rho \widehat \ell_{n,\epsilon}(\rho)
    \widehat \ell_{n,\epsilon}(1-\rho)=2\lambda_n.
  \]
\end{lem}

\begin{proof} Let
  \[
    P_n(t)=\sum_{j=1}^n\binom nj \frac{t^{j-1}}{(j-1)!}
  \]
  and
  \[
    g_n(x)=\begin{cases} P_n(\log x) & \text{if $0<x\leqslant 1$},\\
              0           & \text{if $x>1$}\end{cases}
  \]
  for $n=1,2,\ldots$.

  For $0<\epsilon<1$ we denote
  \[
    p_{n,\epsilon}(x)=\begin{cases}  g_n(x) & \text{if $x>\epsilon$},\\
               0      & \text{if $x\leqslant\epsilon$}\end{cases}
  \]
  and
  \[
    \tau(x)=\begin{cases} \frac{c_0}{\epsilon}
      \exp\left(-1/[1-(\frac{x-1}{\epsilon})^2]\right)
        & \text{if $|x-1|<\epsilon$,}\\
      0 & \text{if $|x-1|\geqslant\epsilon$}
    \end{cases}
  \]
  with $c_0$ given by the identity $\int_0^\infty \tau(x)\,dx=1$.

  We define
  \[
    \ell_{n,\epsilon}(x)=\int_0^\infty p_{n,\epsilon}(xy)\tau(y)\,dy.
  \]
  Then $\ell_{n,\epsilon}(x)$ is a smooth
  function on $\mathbb R$ whose support is contained in the interval
  $(\frac{\epsilon}{1+\epsilon}, \frac{1}{1-\epsilon})$.
  Since
  \begin{equation}\label{eq3.1}
    \widehat \ell_{n,\epsilon}(1-s)
    = \widehat p_{n,\epsilon}(1-s)\widehat\tau (s),
  \end{equation}
  we have
  \begin{multline}\label{eq3.2}
    \widehat \ell_{n,\epsilon}(1-s)
    \widehat \ell_{n,\epsilon}(s)
    -\widehat p_{n,\epsilon}(1-s)\widehat p_{n,\epsilon}(s)
    = \widehat p_{n,\epsilon}(1-s)\widehat p_{n,\epsilon}(s)
    \left\{\widehat\tau(s)\widehat\tau(1-s)-1\right\}.
  \end{multline}

  By integration by parts $n-1$ times for the following second integral we get
  \begin{multline}\label{eq3.3}
    \widehat p_{n,\epsilon}(s)=\int_0^1 g_n(x)x^{s-1}
    -\int_0^\epsilon g_n(x)x^{s-1}\,dx \\
    = 1-(1-\frac{1}{s})^n-P_n(\log\epsilon)\frac{\epsilon^s}{s}
    +P_n^\prime(\log\epsilon)\frac{\epsilon^s}{s^2}+\cdots+
    (-1)^{n-1}P_n^{(n-2)}(\log\epsilon)\frac{\epsilon^s}{s^{n-1}}
    +(-1)^n n \frac{\epsilon^s}{s^n} \\
    = O\left(\frac{1}{|s|}+|\log\epsilon|^{n-1}\frac{\epsilon^{\Re s}}{|s|}\right)
  \end{multline}
  for $0<\Re s<1$ and $|s|\geqslant 1$.

  For $0<\Re s<1$,
  \begin{equation}\label{eq3.4}
    1-\widehat\tau(s) = c_0\int_{-1}^1 e^\frac{1}{t^2-1}
    \left[1-(1+t\epsilon)^{s-1}\right]\,dt \leq c_0\int_{-1}^1
    e^\frac{1}{t^2-1}\left(1+\frac{1}{1-\epsilon}\right)\,dt \ll 1.
  \end{equation}
  By \eqref{eq3.3} and \eqref{eq3.4},
  \begin{multline*}\sum_\rho \widehat p_{n,\epsilon}
    (1-\rho)\widehat p_{n,\epsilon}(\rho)
    \left(\{\widehat\tau(\rho)(\widehat\tau(1-\rho)-1)+(\widehat\tau(\rho)-1)\right\} \\
    \ll\sum_\rho \left(\frac{1}{|\rho|}
    +|\log\epsilon|^{n-1}\frac{\epsilon^{\Re\rho}}{|\rho|}\right)
    \left(\frac{1}{|1-\rho|}+|\log\epsilon|^{n-1}
    \frac{\epsilon^{1-\Re\rho}}{|1-\rho|}\right) \\
    \times\max\left(\Big|\int_{-1}^1 e^\frac{1}{t^2-1}
    \left[1-(1+t\epsilon)^{-\rho}\right]\,dt\Big|, \Big|\int_{-1}^1 e^\frac{1}{t^2-1}
    \left[1-(1+t\epsilon)^{\rho-1}\right]\,dt\Big|\right).
  \end{multline*}

  Similarly as in the proof of \cite[(3.9), p.~284]{Bombieri-Lagarias},
  by the De La Vall\'ee-Poussin zero-free region we have
  \[
    \frac{c}{\log(|\rho|+2)}\leqslant \Re(\rho)\leqslant 1-\frac{c}{\log(|\rho|+2)}
  \]
  for some constant $c>0$.  Thus we have
  \begin{equation}\label{eq3.5}
    \frac{\epsilon^{Re(\rho)}}{\sqrt{|\rho|}}\leqslant
    \max_\rho\epsilon^{c/\log(|\rho|+2)}|\rho|^{-1/2}
    =O\left(e^{-c^\prime\sqrt{\log(1/\epsilon)}}\right)
  \end{equation}
  for some constant $c^\prime>0$.

  From \eqref{eq3.5} we deduce that
  \begin{multline*}
    \left(\frac{1}{|\rho|}+|\log\epsilon|^{n-1}
    \frac{\epsilon^{\Re\rho}}{|\rho|}\right)\left(\frac{1}{|1-\rho|}+|\log\epsilon|^{n-1}
    \frac{\epsilon^{1-\Re\rho}}{|1-\rho|}\right) \\
    =\frac{1}{|\rho(1-\rho)|}\left\{1+|\log\epsilon|^{n-1}(\epsilon^{\Re\rho}
    +\epsilon^{1-\Re\rho})+|\log\epsilon|^{2n-2}\epsilon\right\}\ll |\rho|^{-3/2}.
  \end{multline*}
  It follows that
  \begin{multline*}
    \sum_\rho \widehat p_{n,\epsilon}(1-\rho)\widehat p_{n,\epsilon}(\rho)
    \left(\{\widehat\tau(\rho)(\widehat\tau(1-\rho)-1)+(\widehat\tau(\rho)-1)\right\}  \\
    \ll\sum_{\rho}\frac{1}{|\rho|^\frac{3}{2}}\max\left(
    \Big|\int_{-1}^1 e^\frac{1}{t^2-1}\left[1-(1+t\epsilon)^{-\rho}\right]\,dt\Big|,
    \Big|\int_{-1}^1 e^\frac{1}{t^2-1}\left[1-(1+t\epsilon)^{\rho-1}\right]\,dt\Big|\right),
  \end{multline*}
  where the right side converges uniformly with respect to sufficiently small
  positive $\epsilon$.  Thus, we can interchange the order of taking
  limit $\epsilon\to 0+$ and summing on $\rho$'s to get that
  \[
    \lim_{\epsilon\to 0+} \sum_\rho \widehat p_{n,\epsilon}(1-\rho)
    \widehat p_{n,\epsilon}(\rho)\left\{\widehat\tau(\rho)
    \widehat\tau(1-\rho)-1\right\}=0.
  \]
  Then it follows from \eqref{eq3.2} that
  \begin{equation}\label{eq3.6}
    \lim_{\epsilon\to 0+}\sum_\rho \widehat \ell_{n,\epsilon}(\rho)
    \widehat \ell_{n,\epsilon}(1-\rho)
    = \lim_{\epsilon\to 0+}\sum_\rho \widehat p_{n,\epsilon}(\rho)
    \widehat p_{n,\epsilon}(1-\rho).
  \end{equation}

  By \eqref{eq3.3} we can write
  \begin{multline}\label{eq3.7}
    \widehat g_n(s)\widehat g_n(1-s)-\widehat p_{n,\epsilon}(s)
    \widehat p_{n,\epsilon}(1-s) \\
    = [\widehat g_n(s)-\widehat p_{n,\epsilon}(s)]\widehat g_n(1-s)
    +\widehat p_{n,\epsilon}(s)[\widehat g_n(1-s)-\widehat p_{n,\epsilon}(1-s)] \\
    = O\left(|\log\epsilon|^{n-1}\frac{\epsilon^{\Re s}}{|s|}\right)
    O\left(\frac{1}{|s-1|}\right)
    +O\left(\frac{1}{|s|}+|\log\epsilon|^{n-1}\frac{\epsilon^{\Re s}}{|s|}\right)
    O\left(|\log\epsilon|^{n-1}\frac{\epsilon^{1-\Re s}}{|1-s|}\right) \\
    \ll \frac{1}{|s(1-s)|}\big[|\log\epsilon|^{n-1}\epsilon^{\Re s}
    +(1+|\log\epsilon|^{n-1}\epsilon^{\Re s})
    |\log\epsilon|^{n-1}\epsilon^{1-\Re s}\big].
  \end{multline}
  It follows from \eqref{eq3.7} and \eqref{eq3.5} that
  \begin{equation}\label{eq3.8}
    \lim_{\epsilon\to 0+}\sum_\rho\left[\widehat g_n(\rho)\widehat g_n(1-\rho)
      -\widehat p_{n,\epsilon}(\rho)\widehat p_{n,\epsilon}(1-\rho)\right]=0.
  \end{equation}
  The stated identity follows from \eqref{eq3.6}, \eqref{eq3.8},
  functional equation of $\zeta(s)$, and the identity
  \[
    \left[1-\left(1-\frac{1}{s}\right)^n\right] \cdot \left[1-\left(1-\frac{1}{1-s}\right)^n\right]
    = \left[1-\left(1-\frac{1}{s}\right)^n\right] + \left[1-\left(1-\frac{1}{1-s}\right)^n\right].
  \]

  This completes the proof of the lemma.
\end{proof}

\begin{proof}[Proof of Theorem \ref{thm1.2}]  Let $a(t)=1/t(t-1)$ and
  \[
    \alpha(t)=
    \begin{cases}
      (a_1t+a_2)e^{a(t)} & \text{if $0<t<1$},\\
      0                  & \text{if $t\leqslant 0$ or $1\leqslant t$}
    \end{cases}
  \]
  with $a_1$, $a_2$ being chosen so that $\widehat\alpha(1)=0$
  and $\widehat\alpha(0)=1$.

  If we denote
  \[
    \vartheta(t) = \sum_{n=1}^\infty (-1)^{n-1}\alpha(n t)
    = \sum_{n=1}^\infty \alpha(n t)
    -2\sum_{n=1}^\infty \alpha(n 2t),
  \]
  by the Poisson summation formula
  \[
    \vartheta(t)=\frac{1}{t}\sum_{n\neq 0}^\infty\mathfrak{F}\alpha\left(\frac{n}{t}\right)
    -\frac{1}{t}\sum_{n\neq 0}^\infty\mathfrak{F}\alpha\left(\frac{n}{2t}\right).
  \]
  This implies that $\vartheta(t)$ is of rapid decay when $t\to 0, \infty$.
  It follows that $\widehat\vartheta(s)$ is an entire function.
  Since
  \[
    \sum_{n=1}^\infty \frac{(-1)^{n-1}}{n^s}=(1-2^{1-s})\zeta(s)
  \]
  for $\Re s>0$, by analytic extension we have
  \[
    \widehat\vartheta(s)=(1-2^{1-s})\zeta(s)\widehat\alpha(s)
  \]
  for complex $s$.

  Let
  \begin{equation}\label{eq3.9}
    g_{n,\epsilon}(x)=\ell_{n,\epsilon}(x)-\frac{1}{\widehat\vartheta_1(0)}
    \int_0^\infty\ell_{n,\epsilon}(x/u)\vartheta_1(u)\frac{du}{u}
  \end{equation}
  and
  \[
    h_{n,\epsilon}(x)=\int_0^\infty g_{n,\epsilon}(xy)g_{n,\epsilon}(y)\,dy,
  \]
  where
  \[
    \vartheta_1(u)=
    \begin{cases}
      \vartheta(u) & \text{if $u>\epsilon$},\\
      0            & \text{if $u\leqslant\epsilon$}
    \end{cases}
  \]
  for sufficiently small $\epsilon>0$.

  Since $\widehat\vartheta(0)\neq 0$, we have $\widehat\vartheta_1(0)\neq 0$
  for sufficiently small $\epsilon>0$.  Also, $\widehat\vartheta(\rho)=0$
  for nontrivial zeros $\rho$ of $\zeta(s)$.  Thus, we can write
  \begin{multline*}
    \widehat h_{n,\epsilon}(\rho)
    = \widehat \ell_{n,\epsilon}(\rho) \bigg\{1-\frac{1}{\widehat\vartheta_1(0)}
    \bigg[\widehat\vartheta(\rho)-\int_0^{\epsilon}\vartheta(x) x^{\rho-1}\,dx\bigg]\bigg\}
    \widehat \ell_{n,\epsilon}(1-\rho)\bigg\{1-\frac{1}{\widehat\vartheta_1(0)}
    \bigg[\widehat\vartheta(1-\rho)-\int_0^{\epsilon}\vartheta(x) x^{-\rho}\,dx\bigg]\bigg\} \\
    = \widehat \ell_{n,\epsilon}(\rho)\widehat \ell_{n,\epsilon}(1-\rho)
    \bigg\{1+\frac{1}{\widehat\vartheta_1(0)}
    \int_0^{\epsilon}\vartheta(x) x^{\rho-1}\,dx\bigg\}
    \bigg\{1+\frac{1}{\widehat\vartheta_1(0)}
    \int_0^{\epsilon}\vartheta(x)x^{-\rho}\,dx\bigg\}.
  \end{multline*}
  Hence,
  \begin{multline}\label{eq3.10}
    \widehat h_{n,\epsilon}(\rho)
    -\widehat \ell_{n,\epsilon}(\rho)\widehat \ell_{n,\epsilon}(1-\rho)
    = \widehat \ell_{n,\epsilon}(\rho)\widehat \ell_{n,\epsilon}(1-\rho)
    \frac{1}{\widehat\vartheta_1(0)} \bigg\{\int_0^{\epsilon}\vartheta(x) x^{\rho-1}\,dx \\
    + \int_0^{\epsilon}\vartheta(x)x^{-\rho}\,dx
    +\frac{1}{\widehat\vartheta_1(0)}\int_0^{\epsilon}\vartheta(x) x^{\rho-1}\,dx
    \int_0^{\epsilon}\vartheta(x)x^{-\rho}\,dx\bigg\}
  \end{multline}

  Since both $x\alpha^\prime(x)$ and
  its Fourier transform vanish at $x=0$, by the
  Poisson summation
  \begin{multline*}
    x\vartheta^\prime(x) = \sum_{n=1}^\infty nx\alpha^\prime(nx)
    -2\sum_{n=1}^\infty n2x\alpha^\prime(n2x)
    = \frac{1}{x}\sum_{n\neq 0}^\infty \mathfrak{F}(u\alpha^\prime(u))\left(\frac{n}{x}\right)
    -\frac{1}{x}\sum_{n\neq 0}^\infty \mathfrak{F}(u\alpha^\prime(u))\left(\frac{n}{2x}\right).
  \end{multline*}
  This implies that $\vartheta^\prime(x)$ is of rapid decay when $x\to 0$.
  Since $\vartheta(x)$ is also of rapid decay when $x\to 0$, we have
  \[
    \max\{|\vartheta(x)|, |\vartheta^\prime(x)|\} \ll |x|^n
  \]
  for any positive integer $n$ as $x\to 0+$.  By integration by parts,
  \begin{equation}\label{eq3.11}
    \int_0^{\epsilon}\vartheta(x)x^{-s}\,dx = \frac{\vartheta(\epsilon)}{1-s}
    +\frac{1}{s-1}\int_0^{\epsilon}\vartheta^\prime(x)x^{1-s}\,dx< \frac{c\epsilon}{|s|}\end{equation}
  for $0<\Re s<1$ and $|s|>2$, where $c$ is an absolute constant independent
  of $s$.

  By \eqref{eq3.1}, \eqref{eq3.3} and \eqref{eq3.4} we have
  \begin{equation}\label{eq3.12}
    \widehat\ell_{n,\epsilon}(s)\ll \frac{1}{|s|}+|\log\epsilon|^{n-1}
    \frac{\epsilon^{\Re s}}{|s|}\ll \frac{|\log\epsilon|^{n-1}}{|s|}
  \end{equation}
  for $0<\Re s<1$, where the implied constant depends only on $n$.

  From \eqref{eq3.10}, \eqref{eq3.11} and \eqref{eq3.12} we derive that
  \[
    \sum_{\rho}\left(\widehat h_{n,\epsilon}(\rho)
    -\widehat \ell_{n,\epsilon}(\rho)\widehat \ell_{n,\epsilon}(1-\rho)\right)
    \ll \epsilon|\log\epsilon|^{2n-2}\sum_\rho\frac{1}{|\rho|^3}\to 0
  \]
  as $\epsilon\to 0+$.  By Lemma \ref{lem3.1},
  \[
    \lim_{\epsilon\to 0+}\Delta(h_{n,\epsilon})=2\lambda_n.
  \]

  Note that $g_{n,\epsilon}(t)=0$ for
  $t\not\in (\epsilon^2/(1+\epsilon), 1/(1-\epsilon))$ and
  $\widehat g_{n,\epsilon}(0)=0$ by \eqref{eq3.9}.
  Thus, from now on we choose $\mu_\epsilon$ in the definition
  of the set $S$ to be any fixed finite number large enough so that
  \begin{equation}\label{eq3.13}
    \mu_\epsilon\geqslant(1+\epsilon)/\epsilon^2.
  \end{equation}

  This completes the proof of Theorem \ref{thm1.2}.
\end{proof}

\section{Proof of Theorem 1.3}

For any element $F$ in $L^2(C_S)$,  as $E_S(S(\mathbb{A}_S))$ is dense
in $L^2(C_S)$ by Lemma \ref{lem2.7} there exists a sequence of
elements $f_n\in S(\mathbb{A}_S)$
such that $E_S(f_n)\to F$ in $L^2(C_S)$.
By the definition of the inner product on $L^2(X_S)$, the $f_n$'s
form a Cauchy sequence in $L^2(X_S)$.
Since $L^2(X_S)$ is a complete Hilbert space,
there exists a unique element $f\in L^2(X_S)$
such that $f_n\to f$ in $L^2(X_S)$.
Thus, we define $E_S^{-1}(F) = \{f(\xi x): \xi\in O_S^*\}$.

\begin{lem}\label{lem4.1}  Let $g(t)=t^{-1}g_{n,\epsilon}(t^{-1})$. Then
  \[
    \trace_{E_S(Q_\Lambda^\perp)}(T_h)
    =-\trace_{L^2(C_S)}\left(P_\Lambda E_S\mathfrak{F}_S E_S^{-1}
    P_\frac{1}{\Lambda}V_S(h)E_S\mathfrak{F}_S^t E_S^{-1}\right).
  \]
\end{lem}

\begin{proof}   Let  $F_i, i=1,2,\ldots$ be an orthonormal
  base of $E_S(Q_\Lambda^\perp)$. By Lemma \ref{lem2.4} and Theorem \ref{thm1.1},
  \[
    \trace_{E_S(Q_\Lambda^\perp)}(T_h) = \sum_{i=1}^\infty
    \langle V_S(h)\left(S_\Lambda -E_S  \mathfrak{F}_S^t P_\Lambda\mathfrak{F}_S E_S^{-1}\right)
    F_i, F_i\rangle.
  \]
  Since $F_i\in E_S(Q_\Lambda^\perp)$, we have $\mathfrak{F}_S E_S^{-1}F_i(x)=0$
  for $|x|>\Lambda$. Hence, as two sets we have
  \[
    P_\Lambda\mathfrak{F}_S E_S^{-1}F_i = \mathfrak{F}_S E_S^{-1}F_i.
  \]
  By computations we find that
  \[
    E_S \mathfrak{F}_S^t P_\Lambda\mathfrak{F}_S E_S^{-1}F_i
    = E_S \mathfrak{F}_S^t \mathfrak{F}_S E_S^{-1}F_i=F_i.
  \]
  Thus,
  \[
    \trace_{E_S(Q_\Lambda^\perp)}(T_h)
    =-\sum_{i=1}^\infty\left\langle V_S(h)(1-S_\Lambda) F_i, F_i\right\rangle
    =-\sum_{i=1}^\infty\left\langle V_S(h)P_\frac{1}{\Lambda} F_i, F_i\right\rangle.
  \]

  Since $E_S \mathfrak{F}_S^t P_\Lambda \mathfrak{F}_S E_S^{-1}$ is the projection
  of $L^2(C_S)$ onto $E_S(Q_\Lambda^\perp)$, by Lemma \ref{lem2.6}
  \begin{multline*}
    -\trace_{E_S(Q_\Lambda^\perp)}(T_h)
    = \trace_{L^2(C_S)}
    \left(V_S(h)P_\frac{1}{\Lambda}E_S\mathfrak{F}_S^t P_\Lambda\mathfrak{F}_S  E_S^{-1}\right) \\
    = \trace_{L^2(C_S)}\left(V_S(h)P_\frac{1}{\Lambda}
    E_S\mathfrak{F}_S^t E_S^{-1}P_\Lambda \cdot E_S\mathfrak{F}_S E_S^{-1}\right) \\
    = \trace_{L^2(C_S)}\left(E_S\mathfrak{F}_S E_S^{-1}\cdot V_S(h)P_\frac{1}{\Lambda}
    E_S\mathfrak{F}_S^t E_S^{-1}P_\Lambda\right) \\
    = \trace_{L^2(C_S)}\left(P_\Lambda E_S\mathfrak{F}_S E_S^{-1}
    P_\frac{1}{\Lambda}V_S(h)E_S\mathfrak{F}_S^t E_S^{-1}\right)
  \end{multline*}
  as $V_S(h)^t = V_S(h)$.

  This completes the proof of the lemma.
\end{proof}

\begin{lem}\label{lem4.2}  If we denote
  $\Phi(z,y) = \int_{\mathbb{A}_S}g(uz)\Psi_S(-uy)\,du$, then we can write
  \begin{multline*}
    P_\Lambda E_S\mathfrak{F}_S E_S^{-1}
    P_\frac{1}{\Lambda}V_S(h)E_S\mathfrak{F}_S^t E_S^{-1} F(x) \\
    = P_\Lambda(x)\int_{\mathbb{A}_S,|v|<\frac{1}{\Lambda}}\Psi_S(xv)\,dv
    \int_0^\infty g(vz)\,dz\int_{C_S}\Phi(z,y)
    \sqrt{|xy|}F(y)\,d^\times y.
  \end{multline*}
\end{lem}

\begin{proof}
  Let $F=E_S(f)$ with $f\in S(\mathbb{A}_S)$. We can write
  $E_S\mathfrak{F}_S^t E_S^{-1}F(u) = E_S(\mathfrak{F}_S^t f)(u)$.  Hence,
  \begin{multline*}
    V_S(h)E_S\mathfrak{F}_S^t E_S^{-1} F(v)
    = \int_{C_S}h(v/u)\sqrt{|v/u|}E_S(\mathfrak{F}_S^t f)(u)\,d^\times u \\
    = \int_{C_S}E_S(\mathfrak{F}_S f)(u)\sqrt{|vu|}\,d^\times u\int_0^\infty g(uz)g(vz)\,dz \\
    = \int_{C_S,\frac{|v|(1-\epsilon)}{\mu_\epsilon}
      <|u|< \frac{\mu_\epsilon|v|}{1-\epsilon}} E_S(\mathfrak{F}_S^t f)(u)
    \sqrt{|vu|}\,d^\times u\int_\frac{1-\epsilon}{|v|}^\frac{\mu_\epsilon}{|v|}g(uz)g(vz)\,dz \\
    =\int_\frac{1-\epsilon}{|v|}^\frac{\mu_\epsilon}{|v|}g(vz)\,dz
    \int_{C_S,\frac{|v|(1-\epsilon)}{\mu_\epsilon}
      <|u|<\frac{\mu_\epsilon|v|}{1-\epsilon}}
    E_S(\mathfrak{F}_S^t f)(u)g(uz)\sqrt{|vu|}\,d^\times u \\
    = \int_0^\infty g(vz)\,dz\int_{C_S}g(uz)E_S(\mathfrak{F}_S f)(u)\sqrt{|vu|}\,d^\times u \\
    = \int_0^\infty \sqrt{|v|}g(vz)\,dz\int_{C_S}g(uz)
    \Big(\sum_{\eta\in O_S^*}\mathfrak{F}_S f(\eta u)\Big)|u|\,d^\times u
  \end{multline*}
  where the change of order of integration between $d^\times u$
  and $dz$ is permissible because the double integral is absolutely integrable
  by the choice of $f$.

  Let $\phi(x)=g(|x|)$ if $x\in I_S$ and
  $\phi(x)=0$ if $x\in\mathbb{A}_S-I_S$.  Then $\phi\in S(\mathbb{A}_S)$.
  For each $x\in J_S$, by Lemma \ref{lem2.1} there exists exactly one
  $\xi\in O_S^*$ such that $\xi x\in I_S$.  This implies that
  \[
    g(x) = \sum_{\xi\in O_S^*}\phi(\xi x)
  \]
  for all $x\in C_S$.  Also, we can write
  \[
    \mathfrak{F}_S g(t)=\int_{\mathbb{A}_S}g(|\lambda|)\Psi_S(-\lambda t)\,d\lambda
    = \sum_{\xi\in O_S^*}\int_{\xi I_S}
    g(|\lambda|)\Psi_S(-\lambda t)\,d\lambda.
  \]
  By using the above two identities we can write
  \begin{multline}\label{eq4.1}
    \int_{C_S}g(uz)\Big(\sum_{\eta\in O_S^*}\mathfrak{F}_S f(\eta u)\Big)|u|\,d^\times u\\
    = \int_{C_S}\Big[\sum_{\xi\in O_S^*}\phi(\xi uz)\Big]
    \Big[\sum_{\eta\in O_S^*}\mathfrak{F}_S f(\eta u)\Big]|u|\,d^\times u
    = \int_{C_S}{1\over|z|}\Big[\sum_{\xi\in O_S^*}\mathfrak{F}_S\phi(\xi {u\over z})\Big]
    \Big[\sum_{\eta\in O_S^*}f(\eta u)\Big]|u|\,d^\times u \\
    = \int_{C_S}\Big[\int_{\mathbb{A}_S}g(uz)\Psi_S(-uy)\,du\Big]
    \Big[\sum_{\xi\in O_S^*}f(\xi y)\Big]|y|\,d^\times y
    = \int_{C_S}\Phi(z,y)F(y)\sqrt{|y|}\,d^\times y,
  \end{multline}
  where the right side of the second equality is obtained
  by using the fact that $\mathfrak{F}_S$ is unitary on $L^2(X_S)$;
  see \cite[Lemma 1 b), p.~55]{Connes}.  Thus, we can write
  \[
    V_S(h)E_S\mathfrak{F}_S^t E_S^{-1} F(v)
    = \int_0^\infty g(vz)\,dz\int_{C_S}\Phi(z,y)
    \sqrt{|v/y|}F(y)|y|\,d^\times y.
  \]
  It follows that
  \begin{multline}\label{eq4.2}
    P_\Lambda E_S\mathfrak{F}_S E_S^{-1}
    P_\frac{1}{\Lambda} V_S(h)E_S\mathfrak{F}_S^t E_S^{-1} F(x) \\
    = P_\Lambda(x)\int_{\mathbb{A}_S,|v|<\frac{1}{\Lambda}} \Psi_S(xv)\,dv
    \int_0^\infty g(vz)\,dz\int_{C_S}\Phi(z,y)
    \sqrt{|xy|}F(y)\,d^\times y.
  \end{multline}

  This completes the proof of the lemma.
\end{proof}

\begin{lem}\label{lem4.3}  We have
  \[
    \trace_{E_S(Q_\Lambda^\perp)}(T_h)
    = -\int_{C_S,|x|<\Lambda}|x|\,d^\times x\int_{\mathbb{A}_S,|v|<\frac{1}{\Lambda}}
    \Psi_S(xv)\,dv\int_0^\infty g(vz)\Phi(z,x)\,dz.
  \]
\end{lem}

\begin{proof}  Since $g(vz)=g(|vz|)$, by changing variables
$v\to v(|x|, 1, \cdots, 1)/x$ we can write
  \[
    \int_{\mathbb{A}_S,|v|<\frac{1}{\Lambda}}
    \Psi_S(xv)\,dv\int_0^\infty g(vz)\Phi(z,x)\,dz
    =\int_{\mathbb{A}_S,|v|<\frac{1}{\Lambda}}
    \Psi_S(v(|x|,1,\cdots, 1))\,dv\int_0^\infty g(vz)\Phi(z,x)\,dz
  \]
  and
  \[
    \Phi(z,x)=\int_{\mathbb{A}_S}g(vz)\Psi_S(-v(|x|,1,\cdots, 1))dv=\Phi(z,|x|).
  \]

Let $F=E_S(f)$ with $f\in S(\mathbb{A}_S)$.  Similarly as
  in \eqref{eq2.2} and \eqref{eq2.4}, by integration by parts
  with respect to $v$ we can write
  \begin{multline}\label{eq4.3}
    \int_{\mathbb{A}_S,|v|<\frac{1}{\Lambda}}\Psi_S((|x|,1,\cdots,1)v)\,dv
    \int_0^\infty g(vz)\,dz\int_{C_S}\Phi(z,y)
    \sqrt{|xy|}F(y)\,d^\times y \\
    = 2\sum_{k,l\in\mathbb N_S} \frac{\mu(k)}{k}\int_0^\frac{1}{\Lambda}
    \cos(2\pi |x|v\frac{l}{k})\,dv
    \int_0^\infty g(vz)\,dz\int_{C_S}\Phi(z,y)
    \sqrt{|xy|}F(y)\,d^\times y \\
    = \sum_{k,l\in\mathbb N_S} \frac{\mu(k)}{\pi |x|l}\bigg\{\sin\left(2\pi |x|\frac{l}{k\Lambda}\right)
    \int_0^\infty g\left(\frac{z}{\Lambda}\right)\,dz
    -\int_0^\frac{1}{\Lambda}\sin\left(2\pi |x|v\frac{l}{k}\right)\,dv
    \int_0^\infty zg^\prime(vz)\,dz\bigg\}\\
    \times\int_{C_S}\Phi(z,y)\sqrt{|xy|}F(y)\,d^\times y,
  \end{multline}
  where changing order of differentiation with respect to $v$ and integration
  with respect to $z$ is permissible as $g(vz)=0$
  if $vz\not\in[1-\epsilon,\mu_\epsilon]$.

  By choosing $c=1/8$ in \eqref{eq2.3} we get that
  \[
    \Phi(z,y) \ll_S |yz^7|^{-1/8}.
  \]
  From this inequality we derive that
  \begin{multline}\label{eq4.4}
    \Big|\sin\Big(2\pi |x|\frac{l}{k\Lambda}\Big)\Big|\int_0^\infty
    \left|g\left(\frac{z}{\Lambda}\right)\right|\,dz
    \int_{C_S}|\Phi(z,y)|
    \sqrt{|xy|}|F(y)|\,d^\times y\\
    \ll_S\int_{C_S}|y|^{-1/8}\sqrt{|xy|}|F(y)|\,d^\times y\ll_S\sqrt{|x|}.
  \end{multline}

  Let $0<\nu<1/8$ be a fixed number.  Then
  $|\sin t|\leqslant |\sin t|^{1-\nu} \leqslant |t|^{1-\nu}$ for all real $t$.
  In particular, we have
  \[
    \Big|\sin\Big(2\pi |x|v\frac{l}{k}\Big)\Big|\leqslant \Big|2\pi xv\frac{l}{k}\Big|^{1-\nu}.
  \]

  By changing variables $z\to z/v$, $u\to uv$ and using above inequalities
  for $\Phi(z,y)$ and $\sin(2\pi |x|v\frac{l}{k})$ we derive that
  \begin{multline}\label{eq4.5}
    \int_0^\frac{1}{\Lambda}\Big|\sin\Big(2\pi |x|v\frac{l}{k}\Big)\Big|\,dv
    \int_0^\infty |zg^\prime(vz)|\,dz\int_{C_S}|\Phi(z,y)|
    \sqrt{|xy|}|F(y)|\,d^\times y \\
    =\int_0^\frac{1}{\Lambda}\Big|\sin\Big(2\pi |x|v\frac{l}{k}\Big)\Big|\frac{dv}{v}
    \int_{1-\epsilon}^{\mu_\epsilon}|zg^\prime(z)|\,dz\int_{C_S}|\Phi(z,vy)|
    \sqrt{|xy|}|F(y)|\,d^\times y \\
    \ll_S\int_0^\frac{1}{\Lambda}\Big|2\pi |x|v\frac{l}{k}\Big|^{1-\nu}\frac{dv}{v}
    \int_{1-\epsilon}^{\mu_\epsilon}|zg^\prime(z)|\,dz\int_{C_S}|vy|^{-1/8}
    \sqrt{|xy|}|F(y)|\,d^\times y \\
    \ll_S (|x|l)^{1-\nu}\int_0^\frac{1}{\Lambda}v^{-\nu-1/8}\,dv
    \int_{C_S}|y|^{-1/8}\sqrt{|xy|}|F(y)|\,d^\times y \ll_S \sqrt{|x|}(|x|l)^{1-\nu}.
  \end{multline}

  From \eqref{eq4.4}--\eqref{eq4.5} and $\sum_{l\in\mathbb N_S}l^{-\nu}<\infty$
  we conclude that the series \eqref{eq4.3} converges absolutely and
  \[
    \sum_{k,l\in\mathbb N_S}\frac{\mu(k)}{k}\int_0^\frac{1}{\Lambda}\cos(2\pi |x|v\frac{l}{k})\,dv
    \int_0^\infty g(vz)\,dz\int_{C_S}\Phi(z,y)
    \sqrt{|xy|}F(y)\,d^\times y\ll_S |x|^{-1/2}.
  \]

  The absolute convergence of \eqref{eq4.3}--\eqref{eq4.5} guarantees
  that we can change the order of integration to move the front three terms
  on the right side of \eqref{eq4.2} into $\int_{C_S}\cdots d^\times y$ and get
  \begin{multline*}
    P_\Lambda E_S\mathfrak{F}_S E_S^{-1}P_\frac{1}{\Lambda}V_S(h)E_S\mathfrak{F}_S^t  E_S^{-1} F(x) \\
    =\int_{C_S}P_\Lambda(x)\sqrt{|xy|}F(y)\,d^\times y\int_{\mathbb{A}_S,|v|<\frac{1}{\Lambda}}
    \Psi_S(xv)\,dv\int_0^\infty g(vz)\Phi(z,y)\,dz.
  \end{multline*}
  Since $P_\Lambda E_S\mathfrak{F}_S E_S^{-1} V_S(h) P_\frac{1}{\Lambda} E_S\mathfrak{F}_S^t E_S^{-1}$
  is bounded, this identity holds for all $F\in L^2(C_S)$.

  As $T_h $ is a trace class Hilbert-Schmidt operator
  on $L^2(C_S)$ by Theorem \ref{thm1.1}, it follows from Lemmas
  \ref{lem4.1} and \ref{lem2.5} that
  \[
    \trace_{E_S(Q_\Lambda^\perp)}(T_h)
    =-\int_{C_S,|x|<\Lambda}|x|\,d^\times x\int_{\mathbb{A}_S,|v|<\frac{1}{\Lambda}}
    \Psi_S(xv)\,dv\int_0^\infty g(vz)\Phi(z,x)\,dz.
  \]

  This completes the proof of the lemma.
\end{proof}

\begin{lem}\label{lem4.4}  Let $\Lambda=1$.  Then we can write
  \[
    \trace_{E_S(Q_\Lambda^\perp)}(T_h)
    =-\int_{\mathbb{A}_S,|v|<1}\,dv\int_{C_S,|x|<1}\Psi_S(xv)|x|
    \,d^\times x\int_0^\infty g(vz)\Phi(z,x)\,dz.
  \]
\end{lem}

\begin{proof}  Similarly as in the proof of Lemma \ref{lem4.3},
  by integration by parts with respect to $v$ we derive
  \begin{multline}\label{eq4.6}
    \int_{\mathbb{A}_S,|v|<1}\Psi_S(xv)\,dv
    \int_0^\infty g(vz)\Phi(z,x)\,dz \\
    =2\sum_{k,l\in\mathbb N_S}\frac{\mu(k)}{k}\int_0^1\cos\Big(2\pi |x|v\frac{l}{k}\Big)\,dv
    \int_0^\infty g(vz)\Phi(z,x)\,dz \\
    =\sum_{k,l\in\mathbb N_S}\frac{\mu(k)}{\pi l|x|}\bigg\{\sin\Big(2\pi |x|\frac{l}{k}\Big)
    \int_0^\infty g(z)\Phi(z,x)\,dz-\int_0^1\sin\Big(2\pi |x|v \frac{l}{k}\Big)\,dv
    \int_0^\infty g^\prime(vz)\Phi(z,x)z\,dz\bigg\} \\
    \leqslant \sum_{k,l\in\mathbb N_S}\frac{|\mu(k)|}{\pi l|x|}
    \Big(2\pi |x|\frac{l}{k}\Big)^{1-\nu}\bigg\{\int_{1-\epsilon}^{\mu_\epsilon}|g(z)\Phi(z,x)|\,dz
    +\int_0^1 v^{1-\nu}\,dv\int_0^\infty|g^\prime(vz)\Phi(z,x)|z\,dz\bigg\}.
  \end{multline}

  By choosing $0<c=\nu<1/8$ in \eqref{eq2.3} we deduce that
  \[
    |\Phi(z,x)|=|\Phi(1,x/z)/z|\ll_S|x/z|^{-\nu}/|z|=|x|^{-\nu}|z|^{\nu-1}.
  \]
  It follows that
  \begin{equation}\label{eq4.7}
    \int_{1-\epsilon}^{\mu_\epsilon}|g(z)\Phi(z,x)|\,dz
    \ll_S |x|^{-\nu}
  \end{equation}
  and
  \begin{multline}\label{eq4.8}
    \int_0^1 v^{1-\nu}\,dv\int_0^\infty|g^\prime(vz)\Phi(z,x)|z\,dz
    \ll_S\int_0^1 v^{1-\nu}\,dv\int_0^\infty|g^\prime(vz)||z/x|^\nu \,dz \\
    \ll_S|x|^{-\nu}\int_0^1 v^{-2\nu}\,dv
    \int_{1-\epsilon}^{\mu_\epsilon}|g^\prime(z)||z|^\nu \,dz
    \ll_S |x|^{-\nu}.
  \end{multline}
  By \eqref{eq4.7}--\eqref{eq4.8}, the series \eqref{eq4.6} converges absolutely
  and $\ll_S |x|^{-2\nu}$.

  By changing variables $v\to v(|x|, 1, \cdots, 1)/x$ and by the absolute convergence
  of \eqref{eq4.6} we can change the order of integration between
  $x$ and $v$ as follows to get
  \begin{multline} \label{eq4.9}
    \int_{C_S,|x|<1}|x|\,d^\times x\int_{\mathbb{A}_S,|v|<1}\Psi_S(xv)\,dv
    \int_0^\infty g(|vz|)\Phi(z,x)\,dz \\
    =\int_{I_S,|x|<1}|x|\,d^\times x\int_{\mathbb{A}_S,|v|<1}\Psi_S((|x|,1,\cdots,1)v)\,dv
    \int_0^\infty g(|vz|)\Phi(z,|x|)\,dz \\
    =\int_0^1\,dx \sum_{k,l\in\mathbb N_S}\frac{\mu(k)}{\pi lx}\bigg\{\sin\Big(2\pi x\frac{l}{k}\Big)
    \int_0^\infty g(z)\Phi(z,x)\,dz-\int_0^1\sin\Big(2\pi xv\frac{l}{k}\Big)\,dv
    \int_0^\infty g^\prime(vz)\Phi(z,x)z\,dz\bigg\} \\
    =\sum_{k,l\in\mathbb N_S}\frac{\mu(k)}{\pi l}\bigg\{\int_0^1\sin\Big(2\pi x\frac{l}{k}\Big)
     \,d^\times x \int_0^\infty g(z)\Phi(z,x)\,dz \\
    -\int_0^1 dv \int_0^1\sin\Big(2\pi xv\frac{l}{k}\Big)\,d^\times x
    \int_0^\infty g^\prime(vz)\Phi(z,x)z\,dz\bigg\}.
  \end{multline}

  Also since
  \[
    \int_0^1\cos\Big(2\pi xv\frac{l}{k}\Big)\,dx=\frac{d}{dv}
    \int_0^1 \frac{\sin(2\pi xv\frac{l}{k}) }{ 2\pi x \frac{l}{k}}\,dx,
  \]
  we can write
  \begin{multline} \label{eq4.10}
    \int_{\mathbb{A}_S,|v|<1}dv\int_{C_S,|x|<1}\Psi_S(xv)|x|
    \,d^\times x\int_0^\infty g(vz)\Phi(z,x)\,dz \\
    =\int_{\mathbb{A}_S,|v|<1}dv\int_{C_S,|x|<1}\Psi_S(xv)|x|
    \,d^\times x\int_0^\infty g(vz)\Phi(z,x)\,dz \\
    = 2\sum_{k,l\in\mathbb N_S}\frac{\mu(k)}{k}\int_0^1 dv
    \int_0^1\cos\Big(2\pi xv\frac{l}{k}\Big)\,dx
    \int_0^\infty g(vz)\Phi(z,x)\,dz \\
    = 2\sum_{k,l\in\mathbb N_S}\frac{\mu(k)}{k}\int_0^1 dv\frac{d}{dv}
    \bigg[\int_0^1 \frac{\sin(2\pi xv\frac{l}{k}) }{ 2\pi x \frac{l}{k}}\,dx\bigg]
    \int_0^\infty g(vz)\Phi(z,x)\,dz \\
    =2\sum_{k,l\in\mathbb N_S} \frac{\mu(k)}{k}
    \bigg\{\int_0^1 \frac{\sin(2\pi x \frac{l}{k}) }{ 2\pi x \frac{l}{k}}\,dx
    \int_0^\infty g(z)\Phi(z,x)\,dz \\
    -\int_0^1 dv\int_0^1 \frac{\sin(2\pi xv \frac{l}{k})
    }{ 2\pi x \frac{l}{k}}\,dx
    \int_0^\infty g^\prime(vz)z\Phi(z,x)\,dz\bigg\} \\
    =\sum_{k,l\in\mathbb N_S} \frac{\mu(k)}{\pi l}
    \bigg\{\int_0^1\sin\Big(2\pi x \frac{l}{k}\Big)\,d^\times x
    \int_0^\infty g(z)\Phi(z,x)\,dz \\
    -\int_0^1 dv \int_0^1\sin\Big(2\pi xv \frac{l}{k}\Big) \,d^\times x
    \int_0^\infty g^\prime(vz)\Phi(z,x)z\,dz\bigg\}.
  \end{multline}
  It follows from \eqref{eq4.9} and \eqref{eq4.10} that
  \begin{multline} \label{eq4.11}
    \int_{C_S,|x|<1}|x|\,d^\times x\int_{\mathbb{A}_S,|v|<1}\Psi_S(xv)\,dv
    \int_0^\infty g(|vz|)\Phi(z,x)\,dz \\
    =\int_{\mathbb{A}_S,|v|<1} dv \int_{C_S,|x|<1}\Psi_S(xv)|x|
    \,d^\times x\int_0^\infty g(vz)\Phi(z,x)\,dz.
  \end{multline}
  Then the stated formula follows from Lemma \ref{lem4.3}.

  This completes the proof of the lemma.
\end{proof}

\begin{proof}[Proof of Theorem \ref{thm1.3}]
  The measure difference between $\mathbb{A}_S$ and $J_S$ is negligible
  for a finite set $S$.  By \eqref{eq4.11} and the absolute convergence of \eqref{eq4.6},
  \begin{multline}\label{eq4.12}
    \int_{\mathbb{A}_S,|v|<1} dv \int_{C_S,|x|<1}\Psi_S(xv)|x|
    \,d^\times x\int_0^\infty g(vz)\Phi(z,x)\,dz \\
    =\int_{J_S,|v|<1} dv\int_{C_S,|x|<1}\Psi_S(xv)|x|
    \,d^\times x\int_0^\infty g(|vz|)\Phi(z,x)\,dz \\
    =\sum_{\xi\in O_S^*}
    \int_{\xi I_S, |v|<1} dv\int_{C_S, |x|<1}\Psi_S(xv)|x|\,d^\times x
    \int_0^\infty g(vz)\Phi(z,x)\,dz
  \end{multline}
  converges absolutely.

  Note that $|\xi|=1$ for all $\xi\in O_S^*$.  Because of the absolute
  convergence of \eqref{eq4.12}, for any disjoint decomposition
  $J_S=\cup_{\xi\in O_S^*}\xi I_S$ we have by Lemma \ref{lem4.4}
  \begin{multline}\label{eq4.13}
    \trace_{E_S(Q_\Lambda^\perp)}(T_h)
    =-\int_{\mathbb{A}_S,|v|<1} dv\int_{C_S,|x|<1}\Psi_S(xv)|x|
    \,d^\times x\int_0^\infty g(vz)\Phi(z,x)\,dz \\
    =-\sum_{\xi\in O_S^*}\int_{I_S,|\xi v|<1}d(\xi v)
    \int_{C_S,|x|<1}\Psi_S(x\xi v)|x|\,d^\times x
    \int_0^\infty g(|\xi vz|)\,dz\int_{\mathbb{A}_S}g(|uz|)\Psi_S(-ux)\,du \\
    =-\sum_{\xi\in O_S^*}\int_{I_S,|v|<1} dv
    \int_{C_S,|x|<1}\Psi_S(x\xi v)|x|\,d^\times x
    \int_0^\infty g(|vz|)\,dz\int_{\mathbb{A}_S}g(|uz|)\Psi_S(-ux)\,du
  \end{multline}
  with the sum \eqref{eq4.13} converging absolutely.

  By changing variables in \eqref{eq4.13} first $x\to\xi^{-1}x$
  and then $u\to u\xi$ we get
  \begin{multline}\label{eq4.14}
    \trace_{E_S(Q_\Lambda^\perp)}(T_h) \\
    =-\sum_{\xi\in O_S^*}\int_{I_S,|v|<1} dv
    \int_{C_S,|x|<1}\Psi_S(xv)|x|\,d^\times x
    \int_0^\infty g(|vz|)\,dz\int_{\mathbb{A}_S}g(|uz|)\Psi_S(-ux)\,du,
  \end{multline}
  where \eqref{eq4.14} converges absolutely and sums the same number
  infinitely many times.

  Since the sum \eqref{eq4.14} is finite by Lemma \ref{lem4.1},
  we must have
  \begin{equation}\label{eq4.17}
    \int_{I_S,|v|<1} dv \int_{C_S,|x|<1}\Psi_S(xv)|x|\,d^\times x
    \int_0^\infty g(|vz|)\,dz\int_{\mathbb{A}_S}g(|uz|)\Psi_S(-ux)\,du=0.
  \end{equation}
  From \eqref{eq4.14} and \eqref{eq4.17} we deduce that
  \[
    \trace_{E_S(Q_\Lambda^\perp)}(T_h)=0.
  \]

  This completes the proof of Theorem \ref{thm1.3}.
\end{proof}

\section{Proof of Theorem \ref{thm1.4}}

\begin{lem}\label{lem5.1} $V_S(h)$ is a positive operator on $L^2(C_S)$.
\end{lem}

\begin{proof}  Let $F$ be any element in $L^2(C_S)$
  with compact support.  By definition,
  \[
    V_S(h) F(x)=\int_{C_S} F(\lambda)\sqrt{|x/\lambda|}\,d^\times\lambda
    \int_0^\infty g(|x/\lambda|y)g(y)\,dy.
  \]
  By changing variables $y\to |\lambda|y$ we can write
  \[
    \int_{C_S} V_S(h) F(x) \bar F(x)\,d^\times x
    =\int_{C_S} \bar  F(x)\sqrt{|x|}\,d^\times x
    \int_{C_S}  F(\lambda)\sqrt{|\lambda|}\,d^\times\lambda
    \int_0^\infty g(|x|y)g(|\lambda| y)\,dy.
  \]
  Since the triple integral above is absolute integrable as $F, g$
  are compactly supported, we can change order of integration to derive
  \[
    \int_{C_S} V_S(h)F(x)\bar F(x)\,d^\times x
    =\int_0^\infty \overline{\bigg(\int_{C_S}
      F(x)g(|x|y)\sqrt{|x|}\,d^\times x\bigg)}\bigg(\int_{C_S} F(\lambda)
    g(|\lambda| y)\sqrt{|\lambda|}\,d^\times\lambda\bigg)\,dy\geq 0,
  \]
  where $g$ is a real-valued function.
  Since compactly supported functions are dense in $L^2(C_S)$ and
  $V_S(h)$ is bounded, we have
  \[
    \langle V_S(h)F, F\rangle\geq 0
  \]
  for all $F\in L^2(C_S)$.

  This completes the proof of the lemma.
\end{proof}

\begin{lem}\label{lem5.2} We have
  \[
    \trace_{E_S(Q_\Lambda)}(T_h )
    \geqslant \trace_{E_S(Q_\Lambda)}\{(1-S_\Lambda)T_h\}.
  \]
\end{lem}

\begin{proof}  Let $F_i, i=1,2,\cdots$ be
  an orthonormal base of $E_S(Q_\Lambda)$.
  By Lemma \ref{lem2.4} and Theorem \ref{thm1.1},
  \[
    \trace_{E_S(Q_\Lambda)}(T_h )=\sum_{i=1}^\infty
    \langle V_S(h)\left(S_\Lambda -E_S  \mathfrak{F}_S^t P_\Lambda\mathfrak{F}_S E_S^{-1}\right)
    F_i, F_i\rangle.
  \]
  Since $F_i\in E_S(Q_\Lambda)$, we have $\mathfrak{F}_S E_S^{-1}F_i(x)=0$ for $|x|<\Lambda$.
  This implies that
  \begin{equation}\label{eq5.1}
    P_\Lambda\mathfrak{F}_S E_S^{-1}F_i(x)=0
  \end{equation}
  for all $x$, and hence
  \begin{equation}\label{eq5.2}
    \trace_{E_S(Q_\Lambda)}(T_h )
    =\sum_{i=1}^\infty\langle V_S(h)S_\Lambda F_i, F_i\rangle.
  \end{equation}

  Since $T_h$ is of trace class, so is $(1-S_\Lambda)T_h$
  as $1-S_\Lambda$ is a bounded linear operator on $L^2(C_S)$.
  It follows from Lemma \ref{lem2.4} that the series
  \[
    \sum_{i=1}^\infty \langle (1-S_\Lambda)V_S(h)
    \left(S_\Lambda -E_S  \mathfrak{F}_S^t P_\Lambda\mathfrak{F}_S E_S^{-1}\right)
    F_i, F_i\rangle=\sum_{i=1}^\infty
    \langle V_S(h)S_\Lambda F_i, (1-S_\Lambda) F_i\rangle
  \]
  converges absolutely.  As the right side of \eqref{eq5.2} is also absolutely
  convergent by Lemma \ref{lem2.4} we can write
  \begin{multline*}  \trace_{E_S(Q_\Lambda)}(T_h )
    =\sum_{i=1}^\infty\langle V_S(h)S_\Lambda F_i, S_\Lambda F_i\rangle
    +\sum_{i=1}^\infty \langle V_S(h)S_\Lambda F_i, (1-S_\Lambda) F_i\rangle \\
    =\sum_{i=1}^\infty\langle V_S(h)S_\Lambda F_i, S_\Lambda F_i\rangle
    +\sum_{i=1}^\infty \langle (1-S_\Lambda)T_h F_i, F_i\rangle \\
    =\sum_{i=1}^\infty\langle V_S(h)S_\Lambda F_i, S_\Lambda F_i\rangle
    +\trace_{E_S(Q_\Lambda)}\{(1-S_\Lambda)T_h\}.
  \end{multline*}
  By Lemma \ref{lem5.1}
  \[
    \langle V_S(h)S_\Lambda F_i, S_\Lambda F_i\rangle\geqslant 0
  \]
  for all $i$.  It follows that
  \[
    \trace_{E_S(Q_\Lambda)}(T_h )
    \geqslant \trace_{E_S(Q_\Lambda)}\{(1-S_\Lambda)T_h\}.
  \]

  This completes the proof of the lemma.
\end{proof}

\begin{lem}\label{lem5.3} Let $g(t)=t^{-1}g_{n,\epsilon}(t^{-1})$.  Then
  \begin{multline*}\trace_{E_S(Q_\Lambda)}\{(1-S_\Lambda)T_h\} \\
    =\int_{C_S,\Lambda<|x|}|x|\,d^\times x
    \int_{\mathbb{A}_S,|u|\leqslant \frac{1}{\Lambda}}\Psi_S(ux)\,du\int_0^\infty g(ut)\,dt
    \int_{\mathbb{A}_S,\frac{1}{\Lambda}<|z|}g(zt)\Psi_S(zx)\,dz.\end{multline*}
\end{lem}

\begin{proof}  Since $E_S(1-\mathfrak{F}_S^t P_\Lambda\mathfrak{F}_S)E_S^{-1}$
  is the orthogonal projection of $L^2(C_S)$
  onto $E_S(Q_\Lambda)$, by \eqref{eq5.1}--\eqref{eq5.2} and Lemma \ref{lem2.6}
  \begin{multline*} \trace_{E_S(Q_\Lambda)}((1-S_\Lambda)T_h)
    =\trace_{L^2(C_S)}\left((1-S_\Lambda)
    V_S(h)S_\Lambda E_S(1-\mathfrak{F}_S^t P_\Lambda\mathfrak{F}_S)E_S^{-1}\right) \\
    =\trace_{L^2(C_S)}\{E_S\mathfrak{F}_S E_S^{-1}(1-S_\Lambda)
    V_S(h)S_\Lambda E_S\mathfrak{F}_S^t E_S^{-1}(1-P_\Lambda)\}.
  \end{multline*}

  Let $F=E_S(f)$ with $f\in S(\mathbb{A}_S)$. We have
  \[
    E_S\mathfrak{F}_S^t E_S^{-1}(1-P_\Lambda)F(z)=
    \int_{\mathbb{A}_S}\sqrt{|z/y|}(1-P_\Lambda(y)F(y)\Psi_S(-yz)\,dy.
  \]
  Then
  \begin{multline*}
    V_S(h)S_\Lambda E_S\mathfrak{F}_S^t E_S^{-1}(1-P_\Lambda)F(u) \\
    =\int_{C_S}h(u/z)S_\Lambda(z)\,d^\times z
    \int_{\mathbb{A}_S}\sqrt{|u/y|}(1-P_\Lambda(y))F(y)\Psi_S(-yz)\,dy \\
    =\int_0^\infty g(ut)\,dt \int_{C_S}S_\Lambda(z)g(zt)|z|\,d^\times z
    \int_{\mathbb{A}_S}\sqrt{|u/y|}(1-P_\Lambda(y))F(y)\Psi_S(-yz)\,dy,
  \end{multline*}
  where changing the order of integration in above third line is permissible
  because $g(ut)g(zt)=0$ if $t\not\in |u|^{-1}[1-\epsilon,\mu_\epsilon]$
  or $|z|\not\in|u|[\frac{1-\epsilon}{\mu_\epsilon},\frac{\mu_\epsilon}{1-\epsilon}]$.

  Next, we can write
  \begin{multline*}
    E_S\mathfrak{F}_S E_S^{-1}(1-S_\Lambda)V_S(h)S_\Lambda
    E_S\mathfrak{F}_S^t E_S^{-1}(1-P_\Lambda) F(x) \\
    =\int_{\mathbb{A}_S}\Psi_S(xu)(1-S_\Lambda(u))\,du
    \int_0^\infty g(ut)\,dt \int_{C_S}S_\Lambda(z)g(zt)|z|\,d^\times z
    \int_{\mathbb{A}_S}\sqrt{|x/y|}(1-P_\Lambda(y))F(y)\Psi_S(-yz)\,dy.
  \end{multline*}

  Similarly as in \eqref{eq4.1}, we can write
  \begin{multline*}\int_{C_S}S_\Lambda(z)g(zt)|z|\,d^\times z
    \int_{\mathbb{A}_S}\sqrt{|x/y|}(1-P_\Lambda(y))F(y)\Psi_S(-yz)\,dy \\
    =\int_{C_S}\bigg[\int_{\frac{1}{\Lambda}<|z|}g(zt)
      \Psi_S(zy)\,dz\bigg]\sqrt{|x/y|}(1-P_\Lambda(y))F(y)|y|\,d^\times y.
  \end{multline*}
  It follows that
  \begin{multline*}
    E_S\mathfrak{F}_S E_S^{-1}(1-S_\Lambda)V_S(h)S_\Lambda
    E_S\mathfrak{F}_S^t E_S^{-1}(1-P_\Lambda) F(x) \\
    =\int_{|u|\leqslant \frac{1}{\Lambda}}
    \Psi_S(ux)\,du\int_0^\infty g(ut)\,dt\int_{C_S}\sqrt{|xy|}
    \bigg[\int_{\frac{1}{\Lambda}<|z|}g(zt)\Psi_S(zy)\,dz\bigg](1-P_\Lambda(y))F(y)\,d^\times y.
  \end{multline*}

  An almost identical argument given in \eqref{eq4.3}--\eqref{eq4.5}
  shows that we can move the
  front two terms of the above integral into $\int_{C_S}\cdots d^\times y$ to get
  \begin{multline*}
    E_S\mathfrak{F}_S E_S^{-1}(1-S_\Lambda)V_S(h)S_\Lambda
    E_S\mathfrak{F}_S^t E_S^{-1}(1-P_\Lambda) F(x) \\
    =\int_{C_S}\sqrt{|xy|}\bigg\{\int_{|u|\leqslant\frac{1}{\Lambda}}
    \Psi_S(ux)\,du\int_0^\infty g(ut)\,dt\bigg\}
    \bigg[\int_{\frac{1}{\Lambda}<|z|}g(zt)\Psi_S(zy)\,dz\bigg](1-P_\Lambda(y))F(y)\,d^\times y.
  \end{multline*}

  Since $E_S\mathfrak{F}_S E_S^{-1}(1-S_\Lambda)V_S(h)S_\Lambda
    E_S\mathfrak{F}_S^t E_S^{-1}(1-P_\Lambda)$
  is bounded, the above identity holds for all elements $F$ in $L^2(C_S)$.
  By Lemma \ref{lem2.5},
  \[
    \trace_{E_S(Q_\Lambda)}(1-S_\Lambda)T_h
    =\int_{C_S,\Lambda<|x|}|x|\,d^\times x
    \int_{\mathbb{A}_S,|u|\leqslant \frac{1}{\Lambda}}\Psi_S(ux)\,du\int_0^\infty g(ut)\,dt
    \int_{\mathbb{A}_S,\frac{1}{\Lambda}<|z|}g(zt)\Psi_S(zx)\,dz.
  \]

  This completes the proof of the lemma.
\end{proof}

\begin{proof}[Proof of Theorem \ref{thm1.4}]
  Choosing $\Lambda=1$ in Lemma \ref{lem5.3} we get
  \begin{multline}\label{eq5.3}
    \trace_{E_S(Q_\Lambda)}\{(1-S_\Lambda)T_h\} \\
    =\int_{C_S,1<|x|}|x|\,d^\times x
    \int_{\mathbb{A}_S,|u|<1}\Psi_S(ux)\,du\int_0^\infty g(ut)\,dt
    \int_{\mathbb{A}_S,1<|z|}g(zt)\Psi_S(zx)\,dz,
  \end{multline}
  where we can assume that $1-\epsilon<|ut|<\mu_\epsilon$
  and $1-\epsilon<|zt|<\mu_\epsilon$
  because $g(ut)g(zt)=0$ if $u, z, t$ do not satisfy these two inequalities
  simultaneously.  By these two inequalities we have
  \[
    \max\bigg(\frac{1-\epsilon}{|u|}, \frac{1-\epsilon}{|z|}\bigg)<|t|
    <\min\bigg(\frac{\mu_\epsilon}{|u|}, \frac{\mu_\epsilon}{|z|}\bigg).
  \]
  Since $|u|<1$ and $1<|z|$ by \eqref{eq5.3}, we have
  \[
    \frac{1-\epsilon}{|u|}<|t|<\frac{\mu_\epsilon}{|z|}.
  \]
  This inequality implies that
  \begin{equation}\label{eq5.4}
    1-\epsilon<|t|<\mu_\epsilon, |z|<\frac{\mu_\epsilon}{1-\epsilon},
    \text{ and } \frac{1-\epsilon}{\mu_\epsilon} < |u|.
  \end{equation}
  By using \eqref{eq5.4} we can write
  \begin{multline}\label{eq5.5}
    \trace_{E_S(Q_\Lambda)}\{(1-S_\Lambda)T_h\} \\
    =\int_{C_S, 1<|x|}|x|\,d^\times x
    \int_{\mathbb{A}_S, \frac{1-\epsilon}{\mu_\epsilon}<|u|<1}
    \Psi_S(ux)\,du \int_{1-\epsilon}^{\mu_\epsilon}g(ut)\,dt
    \int_{\mathbb{A}_S, 1<|z|<\frac{\mu_\epsilon}{1-\epsilon}} g(zt)\Psi_S(zx)\,dz.
  \end{multline}

  For $t\in[1-\epsilon, \mu_\epsilon]$, similarly as in \eqref{eq2.3}
  we obtain that
  \begin{equation}\label{eq5.6}
    \int_{\mathbb{A}_S, 1<|z|< \frac{\mu_\epsilon}{1-\epsilon}}
    g(zt)\Psi_S(zx)\,dz\ll_S |x|^{-1}.
  \end{equation}
  Similarly as in \eqref{eq2.2} and \eqref{eq2.4}, by \eqref{eq5.6} and
  integration by parts with respect to variable $u$ we find that
  \begin{multline} \label{eq5.7}
    \int_{\mathbb{A}_S, \frac{1-\epsilon}{\mu_\epsilon}<|u|<1} \Psi_S(ux)\,du
    \int_{1-\epsilon}^{\mu_\epsilon} g(ut)\,dt
    \int_{\mathbb{A}_S, 1<|z|<\frac{\mu_\epsilon}{1-\epsilon}} g(zt)\Psi_S(zx)\,dz \\
    =\sum_{k,l\in N_S} \frac{\mu(k)}{\pi lx} \int_{1-\epsilon}^{\mu_\epsilon}
    \bigg\{g(t)\sin\Big(2\pi x \frac{l}{k}\Big)-t\int_{\frac{1-\epsilon}{\mu_\epsilon}<|u|<1}
    g^\prime(ut)\sin(2\pi ux \frac{l}{k})\,du\bigg\}\,dt\\
    \times\int_{\mathbb{A}_S, 1<|z|< \frac{\mu_\epsilon}{1-\epsilon}} g(zt)\Psi_S(zx)\,dz
    \ll \frac{1}{|x|^2} \sum_{k,l\in N_S} \frac{|\mu(k)|}{l} \ll_S |x|^{-2}.
  \end{multline}

  The above inequality implies that the series
  \begin{multline*}
    \int_{\mathbb{A}_S, \frac{1-\epsilon}{\mu_\epsilon}<|u|<1}\Psi_S(ux)\,du
    \int_{1-\epsilon}^{\mu_\epsilon}g(ut)\,dt
    \int_{\mathbb{A}_S, 1<|z|<\frac{\mu_\epsilon}{1-\epsilon}}g(zt)\Psi_S(zx)\,dz \\
    =\sum_{\gamma\in O_S^*}
    \int_{\gamma I_S, \frac{1-\epsilon}{\mu_\epsilon}<|u|<1}\Psi_S(ux)\,du
    \int_{1-\epsilon}^{\mu_\epsilon}g(ut)\,dt
    \int_{\mathbb{A}_S, 1<|z|<\frac{\mu_\epsilon}{1-\epsilon}}g(zt)\Psi_S(zx)\,dz\ll_S |x|^{-2}
  \end{multline*}
  converges absolutely and uniformly with respect to $|x|>1$.
  By \eqref{eq5.6} and \eqref{eq5.7} we can change the order of
  integration and write \eqref{eq5.5} as
  \begin{multline}\label{eq5.8}
    \trace_{E_S(Q_\Lambda)}\{(1-S_\Lambda)T_h\} \\
    =\int_{\mathbb{A}_S, \frac{1-\epsilon}{\mu_\epsilon}<|u|<1}\,du\int_{C_S, 1<|x|}
    \Psi_S(ux)|x|\,d^\times x \int_{1-\epsilon}^{\mu_\epsilon}g(|ut|)
    [\int_{\mathbb{A}_S, 1<|z|< \frac{\mu_\epsilon}{1-\epsilon}}g(|zt|)\Psi_S(zx)\,dz]\,dt \\
    =\sum_{\gamma\in O_S^*}\int_{I_S, \frac{1-\epsilon}{\mu_\epsilon}<|u|<1}\,du
    \int_{C_S, 1<|x|}\Psi_S(u\gamma x)|x|\,d^\times x
    \int_{1-\epsilon}^{\mu_\epsilon}g(|ut|)\,dt
    \int_{\mathbb{A}_S, 1<|z|<\frac{\mu_\epsilon}{1-\epsilon}}g(|zt|)\Psi_S(zx)\,dz
  \end{multline}
  with the sum \eqref{eq5.8} converging absolutely.

  By changing variables in \eqref{eq5.8} first $x\to x/\gamma$ and
  then $z\to z\gamma$ we deduce that
  \begin{multline}\label{eq5.9}
    \trace_{E_S(Q_\Lambda)}\{(1-S_\Lambda)T_h\} \\
    =\sum_{\gamma\in O_S^*}
    \int_{I_S, \frac{1-\epsilon}{\mu_\epsilon}<|u|<1}\,du\int_{C_S, 1<|x|}
    \Psi_S(ux)|x|\,d^\times x \int_{1-\epsilon}^{\mu_\epsilon}g(|ut|)\,dt
    \int_{\mathbb{A}_S, 1<|z|<\frac{\mu_\epsilon}{1-\epsilon}}g(|zt|)\Psi_S(zx)\,dz,
  \end{multline}
  where \eqref{eq5.9} sums the same number infinitely many times.

  Since the sum \eqref{eq5.9} is finite by Lemma \ref{lem5.3},
  we must have
  \begin{multline}\label{eq5.10}
    \int_{I_S, \frac{1-\epsilon}{\mu_\epsilon}<|u|<1}\,du
    \int_{C_S, 1<|x|} \Psi_S(ux)|x|\,d^\times x \int_{1-\epsilon}^{\mu_\epsilon}g(|ut|)\,dt
    \int_{\mathbb{A}_S, 1<|z|<\frac{\mu_\epsilon}{1-\epsilon}}g(|zt|)\Psi_S(zx)\,dz=0.
  \end{multline}
  Combing \eqref{eq5.9} and \eqref{eq5.10} we get that
  \[
    \trace_{E_S(Q_\Lambda)}\{(1-S_\Lambda)T_h\}=0.
  \]
  By Lemma \ref{lem5.2},
  \[
    \trace_{E_S(Q_\Lambda)}(T_h)\geqslant 0.
  \]

  This completes the proof of Theorem \ref{thm1.4}.
\end{proof}

\section{Proof of Theorem \ref{thm1.5}}

\begin{proof}[Proof of Theorem \ref{thm1.5}]
  By Theorems 1.1--1.4,
  \[
    \Delta(h)=\trace_{E_S(Q_\Lambda^\perp)}(T_h)
    +\trace_{E_S(Q_\Lambda)}(T_h )\geqslant 0.
  \]
  Since
  \[
    h_{n,\epsilon}(x)=\int_0^\infty g_{n,\epsilon}(xy)g_{n,\epsilon}(y)\,dy
    =\int_0^\infty \frac{1}{xy} g_{n,\epsilon}\left(\frac{1}{xy}\right)
    \frac{1}{y}g_{n,\epsilon}\left(\frac{1}{y}\right)\,dy=h(x),
  \]
  we have
  \[
    \Delta(h_{n,\epsilon})=\Delta(h)\geq 0.
  \]
  From Theorem \ref{thm1.2} we deduce that $\lambda_n\geq 0$
  for $n=1,2,\ldots$.  Then the Riemann hypothesis
  \cite[p.~148]{Riemann} follows from Li's criterion
  \cite{Li1} which states that a necessary
  and sufficient condition for the nontrivial zeros of the
  Riemann zeta-function to lie on the critical line is that
  $\lambda_n$ is nonnegative for every positive integer $n$.

  This completes the proof of Theorem \ref{thm1.5}.
\end{proof}


\vskip0.15truein
\noindent{Department of Mathematics, Brigham Young University, Provo, Utah 84602, USA}

\noindent {e-mail: xianjin@mathematics.byu.edu}
\end{document}